\documentclass{amsart}

\usepackage{amsmath,amssymb,amsthm,pinlabel}

\def\co{\colon\thinspace}
\newcommand{\Int}{\mbox{\rm Int}}

\newcommand{\R}{\mathbb{R}}
\newcommand{\Z}{\mathbb{Z}}

\newcommand{\C}{\mathbb{C}}
\newcommand{\LL}{\mathbb{L}}
\newcommand{\tb}{{\tt tb}}
\newcommand{\tS}{{\tt t}_{\Sigma}}
\newcommand{\tD}{{\tt t}_D}
\newcommand{\tT}{{\tt t}_T}

\newcommand{\oL}{\overline{L}}

\newcommand{\rot}{{\tt rot}}
\newcommand{\xist}{\xi_{\mathrm{st}}}
\newcommand{\bandsum}{\#_{\mathrm{Lb}}}

\newtheorem{thm}{Theorem}

\newtheorem{prop}[thm]{Proposition}

\theoremstyle{definition}
\newtheorem*{rem}{Remark}

\newtheorem*{ack}{Acknowledgements}

\begin{document}

\title{Handle moves in contact surgery diagrams}

\author{Fan Ding}
\address{Department of Mathematics, Peking University,
Beijing 100871, P.~R. China}
\email{dingfan@math.pku.edu.cn}
\author{Hansj\"org Geiges}
\address{Mathematisches Institut, Universit\"at zu K\"oln,
Weyertal 86--90, 50931 K\"oln, Germany}
\email{geiges@math.uni-koeln.de}
\date{}

\begin{abstract}
We describe various handle moves in contact surgery diagrams, notably
contact analogues of the Kirby moves. As an application of
these handle moves, we discuss
the respective classifications of long and loose Legendrian knots.
\end{abstract}

\maketitle

\section{Introduction}
In \cite{dige01} and \cite{dige04} we introduced the notion of contact
Dehn surgery, especially contact $(\pm 1)$-surgery, see also
\cite[Section~6.4]{geig08}. The resulting contact surgery diagrams
were studied more carefully in~\cite{dgs04} and since then have led
--- in particular when combined with Heegaard Floer theory ---
to a wealth of applications, see e.g.\
\cite{gls07}, \cite{list04}, \cite{plam06}.
In the present paper we discuss a variety of handle moves in such
contact surgery diagrams: first and second Kirby move, handle
cancellation, and handle moves involving $1$-handles translated
into contact $(+1)$-surgeries. We hope that this will add to
the usefulness of contact surgery diagrams.

The second Kirby move (of which
we shall present three `incarnations') is fundamental to most
of the results in this paper. A further important ingredient
is the Legendrian isotopy between the push-off and the meridian
of a $(-1)$-surgery curve. These two results are proved using
convex surface theory.

Two applications of such handle moves (or the methods used
in their proof) will be discussed: the relation between the classification
of long Legendrian knots in $\R^3$ and their closures in~$S^3$,
and the classification of loose Legendrian knots, i.e.\ knots
with overtwisted complements. The first is in response to
a question by D.~Fuchs and S.~Tabachnikov, the second extends a
result of K.~Dymara.

We assume familiarity with the fundamental notions of
contact topology on the level of~\cite{geig08}. Our contact
structures are understood to be coorientable, i.e.\ they can be written
as $\xi =\ker\alpha$ with $\alpha\wedge d\alpha\neq 0$. With
$\xist$ we denote the standard contact structure $\ker (dz+x\, dy)$
on $\R^3$ as well as the standard contact structure on~$S^3$. Legendrian
knots are always illustrated via their front projections
in the $yz$-plane.
\section{Legendrian band-sum and the second Kirby move}
In this section we formulate the contact analogue of a 2-handle
slide, or what is also called a Kirby move of the second kind.

\begin{figure}[h]
\labellist
\small\hair 2pt
\pinlabel $L_1$ [br] at 89 420
\pinlabel $L_2$ [br] at 446 420
\pinlabel {$L_1\bandsum L_2$} [br] at 106 248
\pinlabel {$L_1\bandsum \oL_2$} [br] at -8 66
\endlabellist
\centering
\includegraphics[scale=0.3]{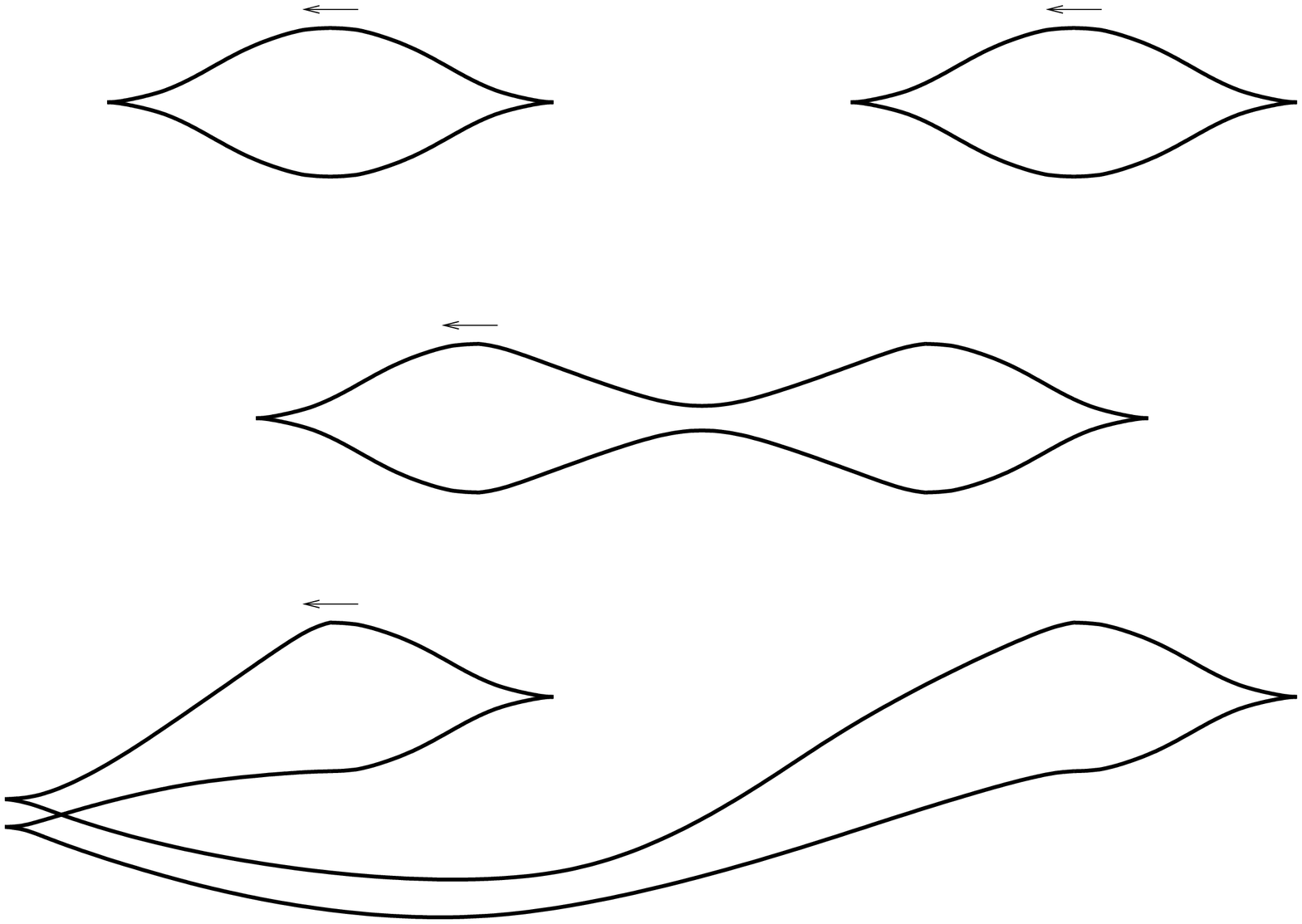}
  \caption{$2$-handle addition and subtraction.}
  \label{figure:2handleslide}
\end{figure}

By a {\em Legendrian band\/} we mean a band in $(\R^3,\xist )$
whose boundary curves consist of an arbitrary Legendrian curve
and its Legendrian push-off in the $z$-direction.
Given two Legendrian knots $L_1$ and $L_2$ in $(\R^3,\xist )$, we
can form their {\em Legendrian band-sum} $L_1\bandsum L_2$.
This means that we construct
the connected sum with the help of a Legendrian band. Here any band
disjoint from the rest of the link may be used. If we think of
$L_1$ and $L_2$ as surgery curves, then --- topologically ---
the replacement of $L_1,L_2$ by $L_1\bandsum L_2,L_2$
corresponds to a handle slide of the $2$-handle corresponding to $L_1$
over that corresponding to~$L_2$. If $L_1$ and $L_2$ are oriented, one can
sensibly speak of {\em handle addition\/} or {\em handle subtraction};
see Figure~\ref{figure:2handleslide}, where we write $\oL_2$ for
$L_2$ with reversed orientation (and the knots are drawn as trivial knots,
but the picture is meant to represent arbitrary knots, which may
also be linked).

\begin{rem}
There is an alternative way for forming the connected sum of two
Legendrian knots in the front projection picture, see
\cite[Figure~2]{etho03}. This could also be used here.
We shall return to this point in Section~\ref{section:bulb}.
\end{rem}

After a topological handle slide, one has to compute the new surgery framing
of $L_1\# L_2$ from the original surgery framings and the linking number
of $L_1$ and $L_2$ (in the case of handles attached to the $4$-ball~$D^4$),
see \cite[p.~142]{gost99}. Here it is obviously necessary to keep
track of the orientation of $L_2$, and to make the distinction
between handle addition and subtraction. For handles attached to
the boundary of an arbitrary $4$-manifold, one can compute the
new framing with the help of the double-strand notation for
framings.

In the contact geometric setting it turns out that, relative to the contact
framing, we do not have to compute new surgery framings at all;
in particular, we may allow rational surgery framings
for~$L_1$, although only integer framings correspond to
actual topological handle slides. (Beware, though, that contact surgery
with rational coefficient different from $1/k$ with $k\in\Z$ is not uniquely
defined, so one can only say that each particular contact $r$-surgery
along $L_1$ corresponds to an $r$-surgery along a
Legendrian isotopic copy of~$L_1$.)
Notice, moreover, that one can make sense of a Legendrian band in an
arbitrary contact manifold, and we shall formulate the following
proposition accordingly. 

\begin{figure}[h]
\labellist
\small\hair 2pt
\pinlabel $L^+$ [bl] at 92 57
\pinlabel $L$ at 20 22
\pinlabel $L$ at 380 22
\pinlabel $L^-$ [b] at 506 69
\endlabellist
\centering
\includegraphics[scale=0.4]{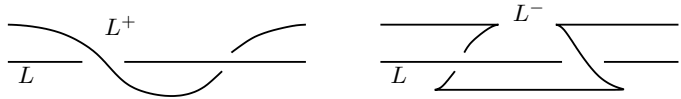}
  \caption{Adding a positive or negative twist.}
  \label{figure:twist}
\end{figure}

\begin{prop}[The second Kirby move]
\label{prop:bandsum}
Let $L_1,L_2$ be two Legendrian knots in a contact $3$-manifold $(M,\xi )$.
Let $L_2^{\pm}$ be the
Legendrian knot obtained as a push-off of $L_2$ with one additional
positive or negative twist, see Figure~\ref{figure:twist}.
Then, in the manifold obtained from $(M,\xi )$ by contact $(\pm 1)$-surgery
along $L_2$, the knot $L_1$ is Legendrian isotopic to
$L_1\bandsum L_2^{\pm}$.
\end{prop}

\begin{proof}
By a {\em Legendrian unknot} $L_0$ we shall mean a Legendrian knot
that can be represented in the front projection of
a suitable Darboux chart by a front with two cusps and no crossings.
The Legendrian band-sum with such a Legendrian unknot is trivial,
provided the band is disjoint from a disc bounded by~$L_0$:
the Legendrian isotopy between $L\bandsum L_0$ and the original Legendrian
knot $L$ is given, in Darboux charts, by a sequence of first and second
Legendrian Reidemeister moves (see~\cite{swia92}
or~\cite{etny05} for these moves). An instance of this phenomenon can
be read off from Figure~\ref{figure:2handleslide}.

Therefore, in order to prove the proposition, it suffices to show
that in the manifold obtained by contact $(\pm 1)$-surgery
along~$L_2$, the knot $L_2^{\pm}$ is a Legendrian unknot as
described.

We are going to write $T$ for any $2$-torus that forms the boundary
of a tubular neighbourhood $\nu L_2$ of~$L_2$. For the computation of
framings we are about to perform, the specific choice of $T$ will
be irrelevant.
A Legendrian push-off $L_2'$ of $L_2$ can be realised as one
of two Legendrian divides on a convex torus~$T$. Along this
Legendrian divide, the contact planes are tangent to~$T$, so
the twisting $\tT (L_2')$ of the contact planes along $L_2'$
relative to the framing induced by $T$ equals~$0$.

The knots $L_2'$ and $L_2^{\pm}$ are topologically isotopic,
so we can compare their framings induced by $\xi$ and $T$, respectively.
As usual, right-handed twists are counted positively. Then the
$\xi$-framing of $L_2^{\pm}$ relative to the $\xi$-framing
of $L_2'$ equals $\pm 1-1$; the $T$-framing of  $L_2^{\pm}$ relative to the
$T$-framing of $L_2'$ equals $\pm 1$. It follows that
\[ \tT (L_2^{\pm})= \tT (L_2')-1=-1.\]

In the surgered manifold, $L_2^{\pm}$ becomes a meridian of the
solid torus we glue in to replace the solid torus $\nu L_2$.
(In particular, $L_2^{\pm}$ is homologically trivial in the
surgered manifold, and we can speak of its Thurston-Bennequin
invariant $\tb (L_2^{\pm})$ there.)
So here the surface framing of $L_2^{\pm}$ given by $T$ coincides
with that given by a meridional disc $D$ with boundary $L_2^{\pm}$.
This implies that
\[  \tb (L_2^{\pm})=\tD (L_2^{\pm})=\tT (L_2^{\pm})=-1\]
in the surgered manifold. Observe that $D$ is disjoint from the
Legendrian band used to connect $L_1$ and $L_2^{\pm}$ in~$M$.

Moreover, $L_2^{\pm}$ can be chosen arbitrarily close to $L_2$.
This guarantees that after the surgery it is a topologically trivial
knot inside a solid torus with tight contact structure. By the
theorem of Eliashberg and Fraser~\cite{elfr98,elfr08}, $L_2^{\pm}$
is a Legendrian unknot in the surgered manifold.
\end{proof}

\begin{rem}
Variants of the second Kirby move will be discussed in
Section~\ref{section:bulb} in the context of the so-called `move~6'
for Legendrian knots, one of the Reidemeister moves in the
presence of $1$-handles.
\end{rem}
\section{Push-off and meridian of $(-1)$-surgery curve}
The aim of the present section is to prove the following proposition, which
will be instrumental in showing how to replace $1$-handles by
contact $(+1)$-surgeries.

\begin{figure}[h]
\labellist
\small\hair 2pt
\pinlabel $-1$ [t] at 30 16
\pinlabel $-1$ [t] at 390 16
\pinlabel $\lambda$ [b] at 250 45
\pinlabel $L$ [t] at 250 16
\pinlabel $L$ [t] at 610 16
\pinlabel $\mu$ [bl] at 575 57
\endlabellist
\centering
\includegraphics[scale=0.4]{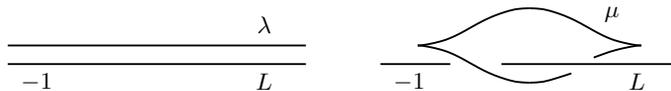}
  \caption{$\lambda$ is Legendrian isotopic to $\mu$ in the surgered manifold.}
  \label{figure:lambda-mu}
\end{figure}

\begin{prop}
\label{prop:lambda-mu}
In the contact manifold obtained from any other contact $3$-manifold
$(M,\xi )$ by contact $(-1)$-surgery along a Legendrian knot $L$, the
Legendrian push-off $\lambda$ of $L$ and a standard Legendrian meridian
$\mu$ (as shown in Figure~\ref{figure:lambda-mu}) are
Legendrian isotopic.
\end{prop}

\begin{proof}
(i) We may assume that $\lambda$ and $\mu$ intersect transversely in
a single point. Choose a tubular neighbourhood $\nu L$ of $L$
such that $\mu ,\lambda \subset \partial (\nu L)$. We claim that
there is a $C^0$-small isotopy of $\partial (\nu L)$ rel $\mu,\lambda$,
turning it into a convex surface, with characteristic foliation
as depicted in Figure~\ref{figure:torus-foliation}; the dashed lines
indicate the dividing curves.

\begin{figure}[h]
\labellist
\small\hair 2pt
\pinlabel $\lambda$ [br] at -3 150
\pinlabel $\mu$ [t] at 185 -3
\endlabellist
\centering
\includegraphics[scale=0.5]{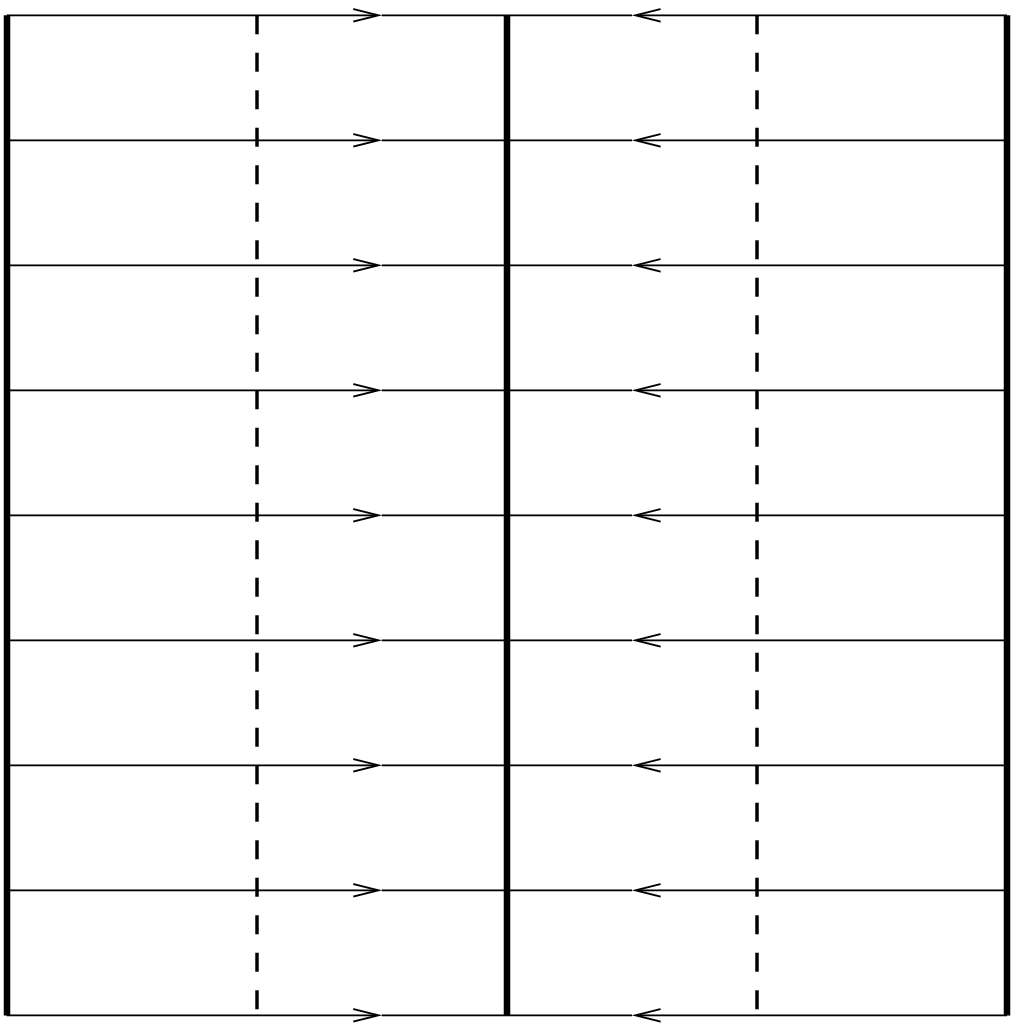}
  \caption{The characteristic foliation on $\partial (\nu L)$.}
  \label{figure:torus-foliation}
\end{figure}

(ii) Assuming this claim for the time being, the proof concludes as follows.
Write $\mu_0:=\ast\times\partial D^2$, $\ast\in S^1$,
and $\lambda_0:=S^1\times\ast$, $\ast\in\partial D^2$,
for meridian and longitude, respectively, of the solid torus $S^1\times D^2$
that we glue in replacing $\nu L$ when we perform surgery along~$L$.
The attaching map for a contact $(-1)$-surgery can be described by
\[ \mu_0\longmapsto \mu-\lambda ,\;\;\; \lambda_0\longmapsto\mu .\]
In other words, we may assume this gluing map to send $\lambda_0$
to a Legendrian ruling curve on $\partial (\nu L)$, and $\lambda_0-\mu_0$
to a Legendrian divide. It follows that the contact structure on
the glued in $S^1\times D^2$, giving us the extension of the
contact structure $\xi$ on $M\setminus\Int (\nu L)$ to that
on the surgered manifold, is given --- up to isotopy rel boundary ---
by the standard model
\[ \cos\theta\, dx -\sin\theta\, dy=0.\]
In this model, the two Legendrian divides $\lambda_0-\mu_0$ as well as
the Legendrian ruling curves $\lambda_0$ are Legendrian isotopic to
the spine $S^1\times\{ 0\}$.

(iii) It remains to prove the claim from (i). Write
$\Sigma :=\partial (\nu L)$, and identify this $2$-torus with
\[ \R /\Z\times \R /\Z \times\{ 0\}\subset
\R /\Z\times \R /\Z \times\R ,\]
with $\mu ,\lambda $ corresponding to the first and second factor
$\R /\Z$, respectively. Write $x,y,z$ for the coordinates corresponding
to these three factors. The twisting $\tS (\lambda )$ of the contact planes
along $\lambda$ relative to the framing induced by $\Sigma$ equals
\[ \tS (\lambda )=0,\]
since $\lambda$ may be realised as one of two Legendrian divides on a
standard convex torus $\Sigma$ with two dividing curves; along such
a Legendrian divide the contact planes coincide with the tangent planes
to~$\Sigma$. Furthermore, with $D$ denoting a meridional disc with
$\partial D=\mu$, we have
\[ \tS(\mu )=\tD (\mu )=\tb (\mu )=-1.\]
It follows that we may perturb $\Sigma$ rel $\mu ,\lambda$
in such a way that the contact structure near the one-point union
of $\mu$ and $\lambda$ looks like
\[ \sin (2\pi x)\, dy+\cos (2\pi x)\, dz =0.\]
So the characteristic foliation looks as in
Figure~\ref{figure:torus-foliation2}~(a), with singular points drawn
in bold face.

\begin{figure}[h]
\labellist
\small\hair 2pt
\pinlabel $y$ [r] at -3 250
\pinlabel $1$ [r] at -3 222
\pinlabel $0$ [tr] at -1 -1
\pinlabel $1$ [b] at 219 -19
\pinlabel $x$ [b] at 250 -19
\pinlabel $y$ [r] at 357 250
\pinlabel $1$ [r] at 357 222
\pinlabel $0$ [tr] at 359 -1
\pinlabel $1$ [b] at 579 -19
\pinlabel $x$ [b] at 610 -19
\pinlabel (a) at 110 255
\pinlabel (b) at 470 255
\endlabellist
\centering
\includegraphics[scale=0.45]{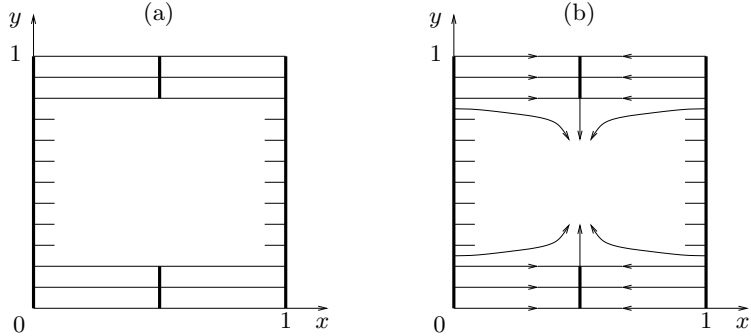}
  \caption{The characteristic foliation on $\Sigma =\partial (\nu L)$
after the first and second perturbation, respectively.}
  \label{figure:torus-foliation2}
\end{figure}

Now perturb $\Sigma$ relative to a small neighbourhood of $\mu$ and $\lambda$
such that $\xi$ looks like
\begin{equation}
\label{eqn:Sigma}
a(x,y)\, dx+\sin (2\pi x)\, dy+\cos (2\pi x)\, dz=0,
\end{equation}
with a smooth function $a(x,y)$ satisfying the following conditions
(for some small $\varepsilon >0$):
\begin{itemize}
\item $a(x,y)=0$ for $y\leq\varepsilon$ or $y\geq 1-\varepsilon$ or $x\leq
\varepsilon$ or $x\geq 1-\varepsilon$,
\item $a(x,y)<0$ for $x$ near $1/2$ and $y\in (\varepsilon ,2\varepsilon )$,
\item $a(x,y)\leq 0$ for $y\in [0,2\varepsilon ]$,
\item analogous conditions near $y=1$ with reversed sign of $a(x,y)$,
\item $\displaystyle{\frac{\partial a}{\partial y}}$ sufficiently small
for (\ref{eqn:Sigma}) to define a contact structure.
\end{itemize}
(The conversion from such a perturbation of $\xi$ to
a perturbation of $\Sigma$ with fixed $\xi$ is a standard argument
using Gray stability, see~\cite[Prop.~4.6.11]{geig08}.)

Then the characteristic foliation on $\Sigma$ is given by the
vector field
\[ X:= \sin (2\pi x)\,\partial_x-a(x,y)\,\partial_y,\]
as indicated in Figure~\ref{figure:torus-foliation2}~(b).
Now perturb $\Sigma$ in the interior (i.e.\ away from
$\mu$ and $\lambda$) to make the characteristic foliation
of Morse-Smale type there. This is possible by the density theorem of
Kupka and Smale, see \cite{peix67} or~\cite{kkn82}. In fact, in the
present situation the vector field $X$ defining the characteristic
foliation near $\mu$ and $\lambda$ is transverse to the boundary of
a disc-shaped region in the interior, and one may appeal directly
to somewhat more elementary results of Peixoto~\cite{peix62} about
Morse-Smale vector fields on surfaces with boundary, with the
vector field being transverse to the boundary.

The resulting surface $\Sigma$ will be convex, cf.\
\cite[Section~4.8]{geig08}. With $\Gamma$ denoting its dividing set,
a result of Kanda~\cite{kand98} allows us to compute the geometric
number of intersection points of $\mu$ and $\lambda$, respectively, with
$\Gamma$ as follows:
\[ -\frac{1}{2}\# (\mu\cap\Gamma )=\tS (\mu )=-1,\;\;\;
-\frac{1}{2}\# (\lambda\cap\Gamma ) =\tS (\lambda )=0.\]
So the dividing set must consist of two curves parallel to $\lambda$.
By a final perturbation of $\Sigma$ (away from
$\mu$ and~$\lambda$) we may achieve the characteristic
foliation depicted in Figure~\ref{figure:torus-foliation},
since that foliation is divided by~$\Gamma$.
\end{proof}
\section{Topologically trivial long Legendrian knots}
A {\em long Legendrian knot\/} is a Legendrian embedding $\gamma\co
\R\rightarrow (\R^3,\xist )$ with $\gamma (t)=(0,t,0)$ for $|t|$ large.
By an isotopy of long knots we mean an isotopy that is stationary
outside a compact subset of~$\R$. A long knot is {\em topologically trivial\/}
if it is isotopic to the embedding $t\mapsto (0,t,0)$. We write
$K=\gamma (S^1)$ for the unparametrised knot; we also use that notation
if the particular parametrisation (except its orientation) is irrelevant.

\begin{figure}[h]
\labellist
\small\hair 2pt
\pinlabel $K$ [t] at 140 -5
\pinlabel $\hat{K}$ [t] at 504 -5
\endlabellist
\centering
\includegraphics[scale=0.45]{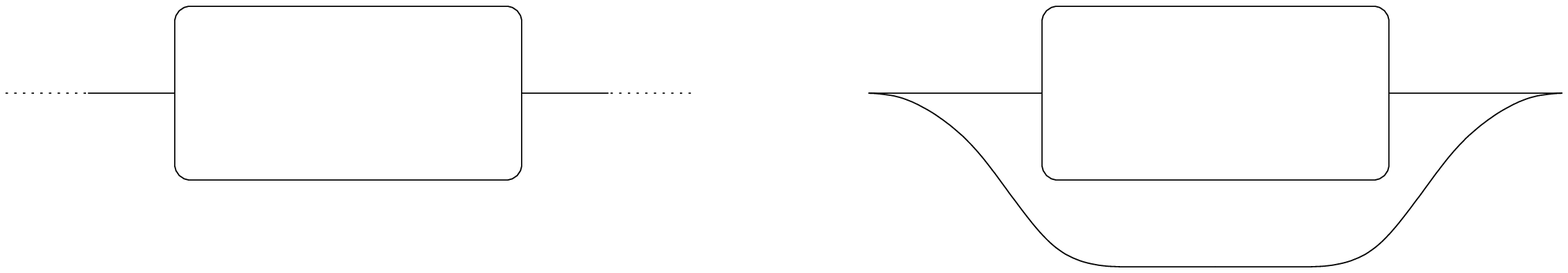}
  \caption{A long knot $K$  and its completion $\hat{K}$.}
  \label{figure:completion}
\end{figure}

The classical invariants $\tb$ (Thurston-Bennequin invariant)
and $\rot$ (rotation number) of a long Legendrian knot are
computed from its front projection in the same way as the invariants
of an honest (`short') Legendrian knot. Any long Legendrian knot $K$
can be completed to a short Legendrian knot $\hat{K}$ as shown in
Figure~\ref{figure:completion}. Then $\tb (\hat{K})=\tb (K)-1$
and $\rot (\hat{K})=\rot (K)$.

Alternatively, one can obtain the completion $\hat{K}$ by regarding
$K$ as a Legendrian knot in $(S^3,\xist )$ passing through
the point at infinity; by a Legendrian isotopy one can then move
$\hat{K}$ into $\R^3\subset S^3$. This is
feasible because the complement of a point in $(S^3,\xist )$ is
contactomorphic to $(\R^3,\xist )$, see \cite[Prop.~2.1.8]{geig08}.

In \cite[Section~2.5]{futa97} the following result is mentioned in passing.

\begin{prop}
\label{prop:trivial-long}
Topologically trivial long Legendrian knots are classified, up to
Legendrian isotopy, by the two classical invariants.
\end{prop}

\begin{proof}
A very short proof of the corresponding result for topologically
trivial Legendrian knots in $(S^3,\xist )$, originally due
to Eliashberg and Fraser~\cite{elfr98,elfr08}, can be found in
\cite[Section~3.5]{etho01}. The argument given there for showing
that any topologically trivial Legendrian knot with
$\tb <-1$ destabilises carries over to long knots
with $\tb <0$ without change. Only the proof
that there is a unique topologically trivial long Legendrian knot with
$\tb =0$ (the maximal possible value of $\tb$ in that case)
requires a  minor modification.

\begin{figure}[h]
\centering
\includegraphics[scale=0.7]{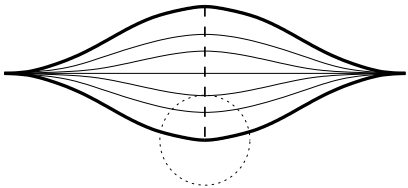}
  \caption{The characteristic foliation on the disc bounded by
the trivial Legendrian knot.}
  \label{figure:saucer-foliation}
\end{figure}

Thus, let $K$ be a topologically trivial {\em long\/} Legendrian knot
in $(\R^3,\xist )$ with $\tb (K)=0$. We regard its completion $\hat{K}$ as
a Legendrian knot in $(S^3,\xist )$ with $\tb (\hat{K}) =-1$.
Let $D$ be a disc bounded
by $\hat{K}$, which we may assume to be convex. Then the dividing set
$\Gamma_D$ consists of a single arc, for otherwise $\hat{K}$ could
be destabilised, contradicting the maximality of~$\tb$. After
a $C^0$-small perturbation of~$D$, the characteristic foliation $D_{\xist}$
may be assumed to look as in Figure~\ref{figure:saucer-foliation},
with two elliptic singular points along~$\hat{K}$
(cf.\ \cite[Prop.~4.8.11]{geig08}).

Given two such knots $K,K'$, which coincide outside
a compact set, we may assume the corresponding discs $D,D'$
to coincide inside a small ball $B^0$ as indicated by a dashed circle
in Figure~\ref{figure:saucer-foliation}.
Notice that one may assume
the characteristic foliation $(\partial B^0)_{\xist}$
to be the standard one with two singular elliptic points
only and no closed orbits.
Since the characteristic foliations $D_{\xist}$ and $D'_{\xist}$ determine
the germ of the contact structure, we find balls $B,B'$ containing
$B^0$ and $D,D'$, respectively, and a contactomorphism
$(B,\xist )\rightarrow (B', \xist )$ fixed on~$B^0$ and
sending $\hat{K}$ to $\hat{K}'$ (and $D$ to~$D'$). By
Eliashberg's classification of tight contact structures on the $3$-ball
\cite[Thm.~2.1.3]{elia92}, this contactomorphism extends to
a contactomorphism of $(S^3,\xist )$. Finally, by the parametric
version of Eliashberg's result and the fact that
the relative diffeomorphism group $\mbox{Diff}(D^3,S^2)$ is
connected~\cite{cerf68}, this contactomorphism of $(S^3,\xist )$
is contact isotopic to the identity relative to~$B^0$. This
translates into a Legendrian isotopy between the long knots
$K$ and~$K'$.
\end{proof}

\begin{rem}
There is a slightly roundabout but more elementary argument for the
final step in the preceding proof. We may choose a smaller ball
$B^1$ inside $B^0$ sharing all its properties, and contained
in the neighbourhood $\nu L_0$ of a standard Legendrian unknot
$L_0$ inside~$B^0$. We then have to deal with a contactomorphism of
$(S^3,\xist )$ fixed on $\nu L_0$. A result of Giroux~\cite{giro01}
can then be used to show that such a contactomorphism is contact isotopic
rel $\nu L_0$ (and in particular rel~$B^1$) to the identity, see
\cite[p.~161/2]{dige07}.
\end{rem}
\section{Cancellation of $1$- and $2$-handles}
\label{section:cancellation}
\begin{figure}[h]
\centering
\includegraphics[scale=0.3]{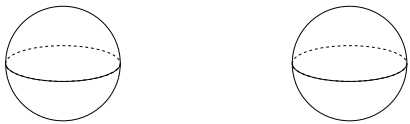}
  \caption{$S^1\times S^2$ with its standard tight contact structure.}
  \label{figure:1-handle}
\end{figure}

\begin{figure}[h]
\labellist
\small\hair 2pt
\pinlabel $-1$ [br] at 385 94
\endlabellist
\centering
\includegraphics[scale=0.3]{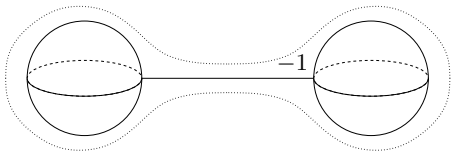}
  \caption{Cancelling handle pair.}
  \label{figure:1-2-handle}
\end{figure}

As described in \cite[Section~2]{gomp98} or \cite[Section 11.1]{gost99},
the diagram in Figure~\ref{figure:1-handle}, depicting a
$1$-handle attached to $D^4$ along $S^3=\partial D^4$, represents
the manifold $S^1\times S^2$ with its standard Stein fillable
contact structure. Since contact $(-1)$-surgery along
a Legendrian knot preserves Stein fillability, the diagram in
Figure~\ref{figure:1-2-handle} represents a Stein fillable contact
structure on $S^3$, so this is a diagram for $(S^3,\xist )$.
A slightly stronger formulation of this observation is that the $3$-ball
bounded by a $2$-sphere $S$ around this configuration (indicated by a dashed
line) carries the tight contact structure uniquely determined,
up to isotopy, by the characteristic foliation~$S_{\xist}$.

We now want to tie up this observation with Propositions
\ref{prop:lambda-mu} and \ref{prop:trivial-long} in order to
prove the following result, which tells us how to replace
$1$-handles by contact $(+1)$-surgeries along trivial Legendrian
knots. This corresponds to the `dotted circle' notation in
Kirby diagrams~\cite[p.~167/8]{gost99}; see also
\cite[p.~167/8]{ozst04} for a related discussion. (It is
interesting to observe the analogy on the level of handlebodies:
a dotted circle
indicates a $2$-handle being `dug out' from~$D^4$, while a 
contact $(+1)$-surgery corresponds to attaching a $2$-handle to
the concave end of a symplectic cobordism.)

The first diagram described in the next theorem is
said to be in `standard form'; by~\cite[Thm.~2.2]{gomp98} any diagram
involving $1$-handles can be written in that way.

\begin{figure}[h]
\labellist
\small\hair 2pt
\pinlabel Legendrian [t] at 324 440
\pinlabel tangle [t] at 324 420
\pinlabel Legendrian [t] at 324 98 
\pinlabel tangle [t] at 324 78
\pinlabel $+1$ [bl] at 279 268
\pinlabel $+1$ [bl] at 279 205
\pinlabel {\LARGE $\cong$} [t] at 326 318
\endlabellist
\centering
\includegraphics[scale=0.5]{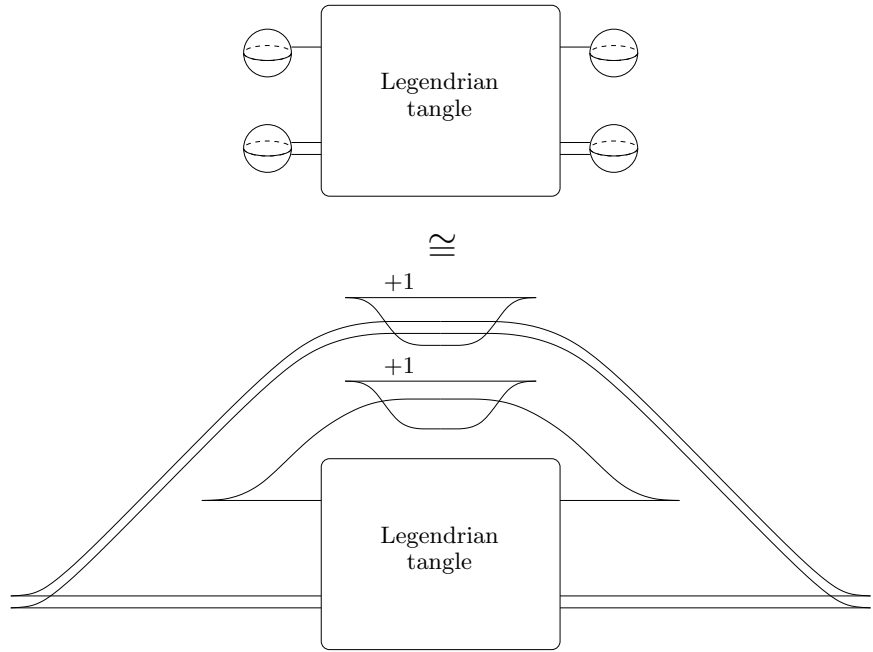}
  \caption{Replacing $1$-handles by $(+1)$-surgeries.}
  \label{figure:link-diagram}
\end{figure}

\begin{thm}
\label{thm:1-to-2}
The two surgery diagrams in Figure~\ref{figure:link-diagram}
represent contactomorphic contact manifolds. Here it is understood
that the surgery coefficients (relative to the contact framing)
of the curves making up the Legendrian tangle remain unchanged.
\end{thm}

\begin{figure}[h]
\labellist
\small\hair 2pt
\pinlabel $-1$ at 142 390
\pinlabel $-1$ at 142 244
\pinlabel $+1$ at 154 357
\pinlabel $+1$ at 154 211
\pinlabel $+1$ at 154 72 
\endlabellist
\centering
\includegraphics[scale=0.45]{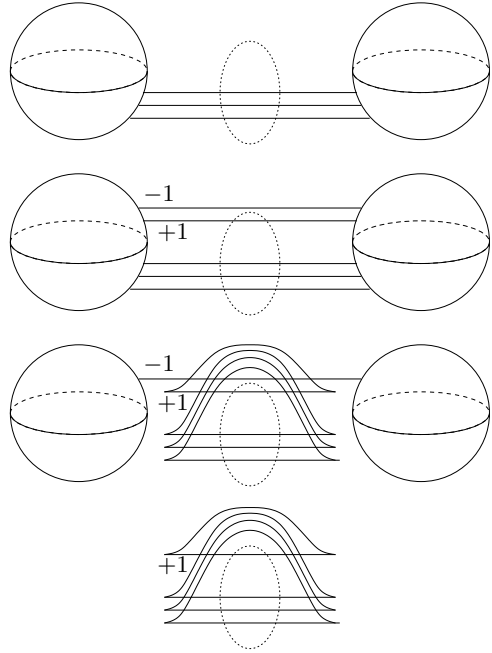}
  \caption{Equivalent surgery diagrams.}
  \label{figure:link-diagram2}
\end{figure}

\begin{proof}
By the cancellation lemma \cite[Prop.~8]{dige01}, cf.\
\cite[Prop.~6.4.5]{geig08}, the first two diagrams shown in
Figure \ref{figure:link-diagram2} are equivalent. The Legendrian
knots in the second diagram are Legendrian isotopic, relative to
the ball indicated by a dashed ellipse, to the knots in
the third diagram by Propositions \ref{prop:lambda-mu}
and~\ref{prop:trivial-long}.
Finally, the fourth diagram in Figure \ref{figure:link-diagram2}
is obtained from the third by cancelling a $1$-handle with a $2$-handle
as described above. The theorem follows by applying this procedure to
each $1$-handle.
\end{proof}

\begin{figure}[h]
\labellist
\small\hair 2pt
\pinlabel $L$ [l] at 306 237
\pinlabel $L_+$ [b] at 511 204
\pinlabel $L_-$ [t] at 511 194
\pinlabel $+1$ [bl] at 128 132
\pinlabel $+1$ [bl] at 128 44
\pinlabel $+1$ [bl] at 525 132
\pinlabel $+1$ [bl] at 525 44
\pinlabel $L_-$ [b] at 196 102
\pinlabel $-1$ [b] at 604 102
\pinlabel $-1$ [b] at 301 14
\pinlabel {$L\bandsum L_-^-$} [l] at 701 147
\pinlabel $\cong$ at 379 199
\pinlabel $\cong$ at 379 19
\endlabellist
\centering
\includegraphics[scale=0.45]{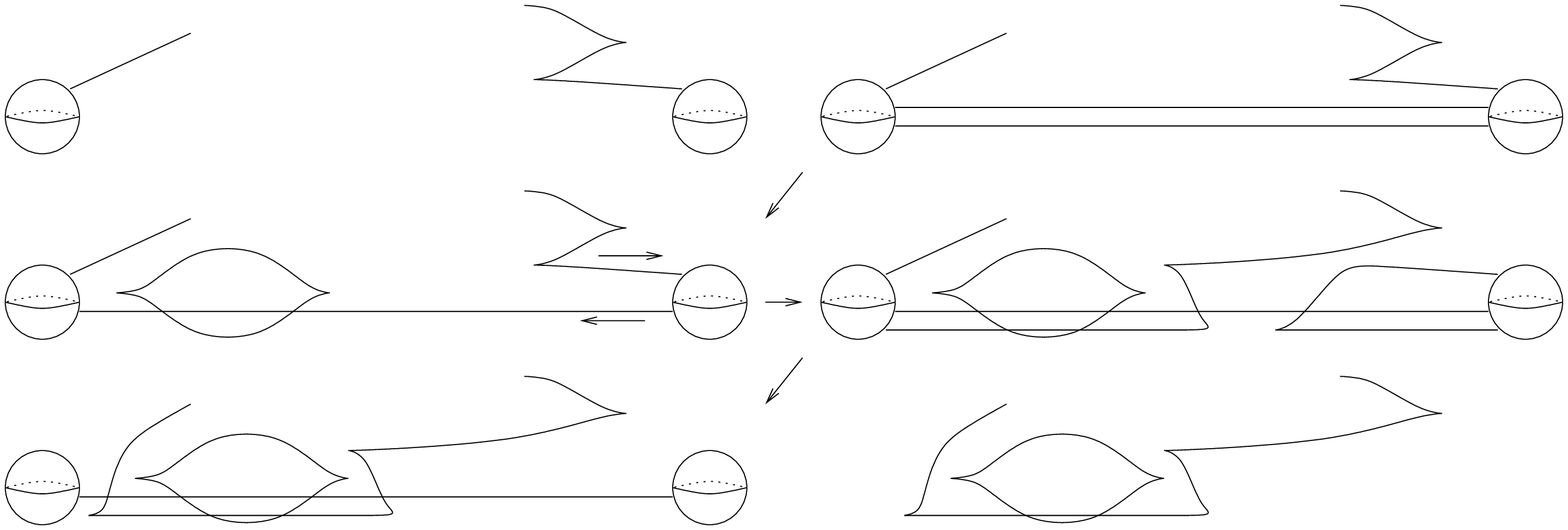}
  \caption{Replacing $1$-handles by $(+1)$-surgeries
in stabilised case.}
  \label{figure:link-diagram-stab}
\end{figure}

When we are dealing with a link diagram where each strand passing
over a $1$-handle is stabilised, one can give a more direct
proof of Theorem~\ref{thm:1-to-2}, using only handle slides
and the cancellation lemma from~\cite{dige04},
cf.~\cite[Prop.~6.4.5]{geig08}, but none of the
awkward arguments involving relative isotopies. This is illustrated
in Figure~\ref{figure:link-diagram-stab} for a single knot~$L$;
the argument for several stabilised strands passing over
a given $1$-handle is analogous.
(In this sequence of
pictures, `$\cong$' indicates a contactomorphism; an arrow,
a Legendrian isotopy.)
First we introduce a
pair of cancelling contact $(\pm 1)$-surgeries along
the Legendrian knots $L_{\pm}$. Then apply
Proposition~\ref{prop:lambda-mu} to turn $L_+$ 
into a meridian of $L_-$. Now form $L\bandsum L_-^-$ as in
Proposition~\ref{prop:bandsum}. The obvious Legendrian isotopy that
gives the fifth picture in the sequence of
Figure~\ref{figure:link-diagram-stab} is a combination of
what are called `move~4' and `move~5'
in~\cite{gomp98}; see Section~\ref{section:4and5} below.
To get the final picture, apply the cancellation
from Figure~\ref{figure:1-2-handle}.

The particular positioning of the meridional $(+1)$-surgery curve
is irrelevant, see the following section. The requirement
that $L$ be stabilised is not restrictive; this will be discussed
in Section~\ref{section:bulb}.

\section{Sliding $1$-handles}
\label{section:1handles}
\begin{figure}[h]
\centering
\includegraphics[scale=0.4]{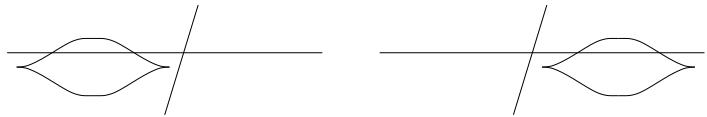}
  \caption{Isotopic Legendrian meridians I.}
  \label{figure:meridian1}
\end{figure}

\begin{figure}[h]
\centering
\includegraphics[scale=0.4]{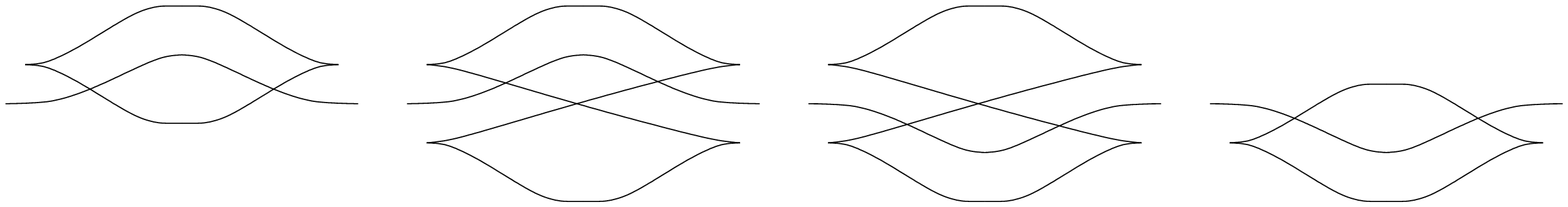}
  \caption{Isotopic Legendrian meridians II.}
  \label{figure:meridian2}
\end{figure}

\begin{figure}[h]
\centering
\includegraphics[scale=0.3]{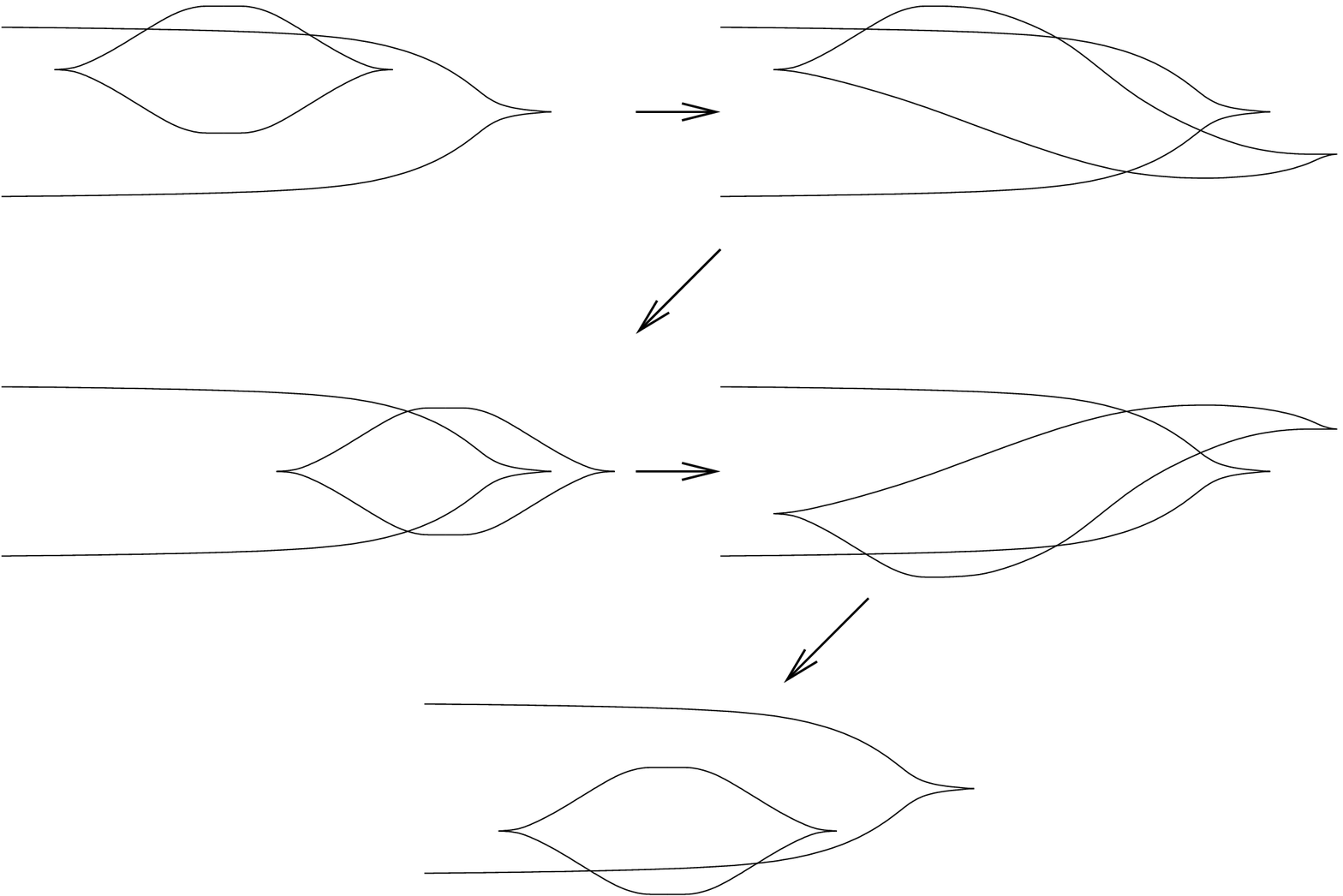}
  \caption{Isotopic Legendrian meridians III.}
  \label{figure:meridian3}
\end{figure}

\begin{figure}[h]
\labellist
\small\hair 2pt
\pinlabel $+1$ [bl] at 148 218
\pinlabel $+1$ [bl] at 452 218
\pinlabel $+1$ [bl] at 148 119
\pinlabel $+1$ [bl] at 452 119
\endlabellist
\centering
\includegraphics[scale=0.45]{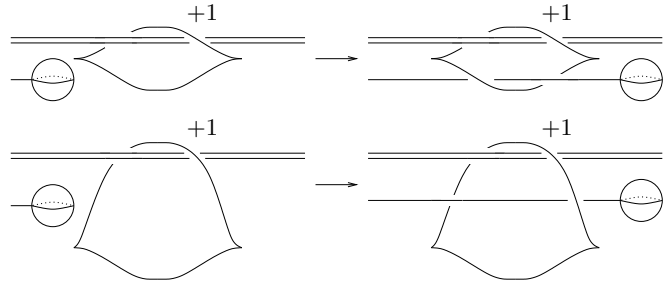}
  \caption{$1$-handle slides of either orientation.}
  \label{figure:1handleslide}
\end{figure}

An easy way to realise $1$-handle slides in a Kirby diagram is to
represent the $1$-handle over which one wants to slide
another $1$-handle by a dotted circle, see \cite[Figure~5.39]{gost99}.
In the Legendrian setting, at first sight there appears to be
less flexibility in performing such moves, because it may not be possible
to move one of the attaching balls of the sliding handle through
the dotted circle by a Legendrian isotopy. However, one can in fact
ensure complete flexibility by moving the dotted circle, which in the
contact setting corresponds to a $(+1)$-surgery along a Legendrian
meridian. Figures \ref{figure:meridian1} to \ref{figure:meridian3}
show Legendrian isotopies of such a meridian;
Figure~\ref{figure:1handleslide} illustrates that a $1$-handle slide of
either orientation can be realised by a contact isotopy.

\section{Moves 4 and 5}
\label{section:4and5}
\begin{figure}[h]
\centering
\includegraphics[scale=0.4]{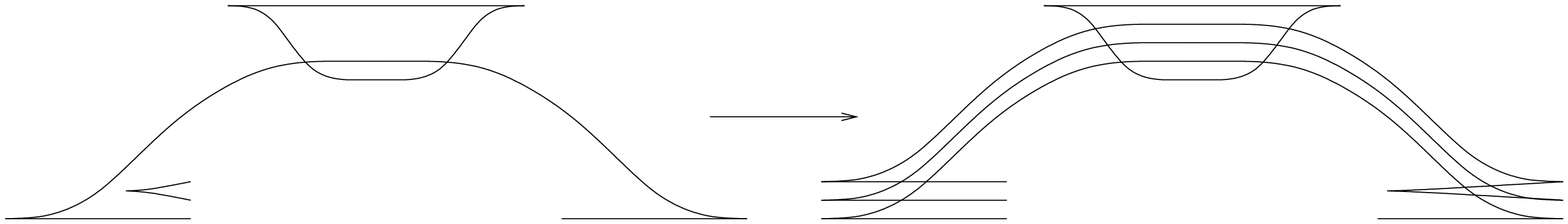}
  \caption{Move 4.}
  \label{figure:move4}
\end{figure}

\begin{figure}
\centering
\includegraphics[scale=0.4]{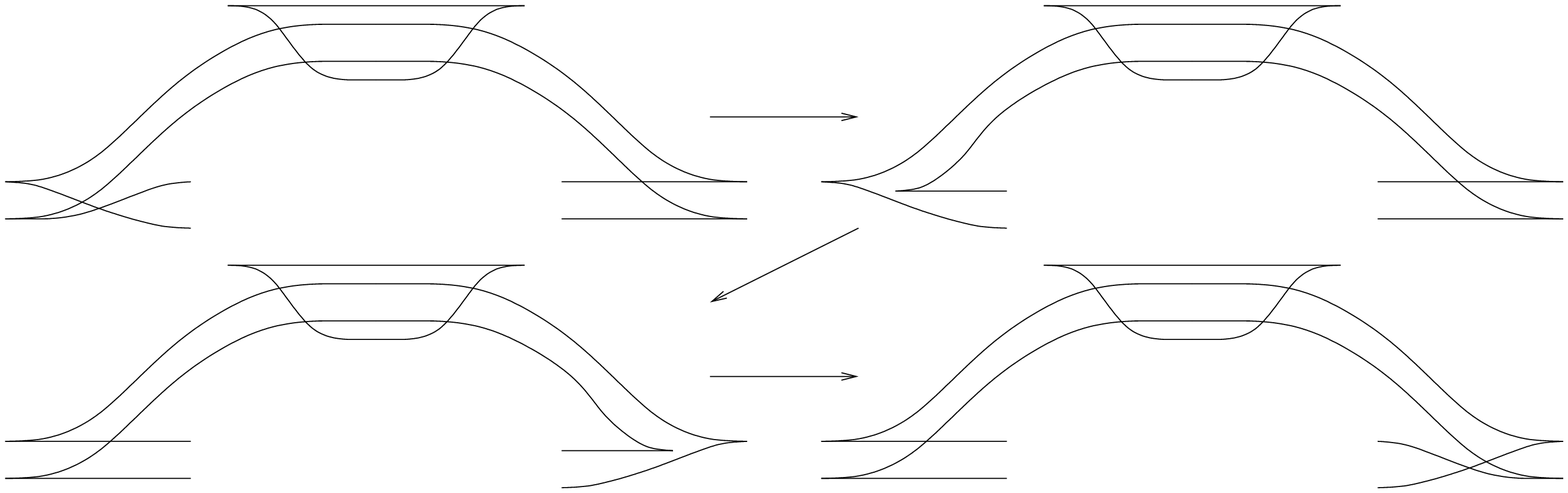}
  \caption{Move 5.}
  \label{figure:move5}
\end{figure}

In addition to the three usual Legendrian Reidemeister moves,
in the presence of $1$-handles there are three further moves
to consider, see~\cite[Thm.~2.2]{gomp98}. We now describe these
moves when the $1$-handle has been replaced by a contact $(+1)$-surgery
along a Legendrian unknot.
Figures \ref{figure:move4} and \ref{figure:move5} give the
simple analogue of move 4 and~5, respectively; move~6 will be
discussed in Section~\ref{section:bulb}. Observe that moves 4 and 5 are
Legendrian isotopies in the complement of the unknot, so the surgery
plays no role here. The isotopy between the first and the second line
in Figure~\ref{figure:move5} is given by two Legendrian
Reidemeister moves of the
second kind (involving only the two strands passing through the unknot)
and a couple of Reidemeister moves of the third kind (involving in addition
the unknot).

\section{Light bulb trick, move 6, and Rolfsen twist}
\label{section:bulb}
\begin{figure}[h]
\centering
\includegraphics[scale=0.45]{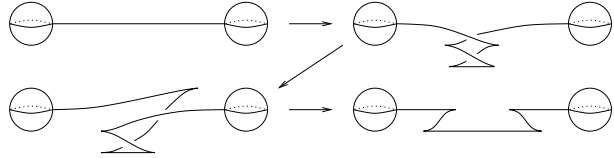}
  \caption{The light bulb trick.}
  \label{figure:lightbulb}
\end{figure}

The requirement at the end of Section~\ref{section:cancellation} that the
knot passing over the $1$-handle be stabilised turns out not
to be any restriction at all. This is seen by a move that
is essentially the light bulb trick~\cite[p.~257]{rolf76},
see Figure~\ref{figure:lightbulb}.
In the contact geometric setting, this move has been
observed earlier by Gompf~\cite[Figure~18]{gomp98}.
The move from the first to the second line in Figure~\ref{figure:lightbulb}
is `move 6' in~\cite{gomp98}.

\begin{rem}
Write $S_{\pm}L$ for the positive or negative stabilisation, respectively,
of an oriented Legendrian knot~$L$, corresponding to adding a down
or up zigzag to the front projection of~$L$. Denote by
\[ \xist =\ker (z\, d\theta +x\, dy-y\, dx).\]
the standard tight (and Stein fillable) contact structure
on $S^1\times S^2\subset S^1\times \R^3$.
Then the light bulb trick says that any Legendrian knot $L$ in
$(S^1\times S^2,\xist )$ that crosses some sphere $\{ \theta_0\}\times
S^2$ exactly once is Legendrian isotopic to its double
stabilisation $S_+S_-L$. Observe that $\xist$ is trivial as a $2$-plane bundle,
with trivialisation given by the vector fields $x\partial_{\theta}-
z\partial_y+y\partial_z$ and $y\partial_{\theta}+z\partial_x-x\partial_z$.

Contrast this with \cite[Thm.~4.1.1]{tche03}, which says that with
$\xist$ replaced by any contact structure of non-zero Euler class,
the knots $L$ and $S_+S_-L$ are not Legendrian isotopic, in spite
of being Legendrian regularly homotopic and framed isotopic; they
can be distinguished by a Vassiliev invariant.
\end{rem}

\begin{figure}[h]
\labellist
\small\hair 2pt
\pinlabel $L_0$ [bl] at 127 120
\pinlabel $L$ [r] at 37 34
\pinlabel {$L\# L_0^+$} [r] at 252 34
\endlabellist
\centering
\includegraphics[scale=0.45]{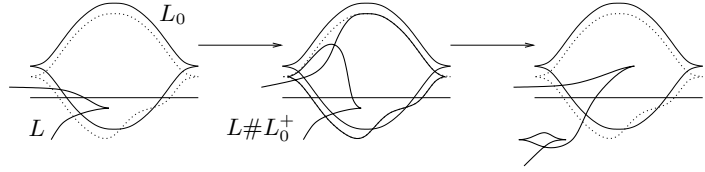}
  \caption{Move 6 via connected sum.}
  \label{figure:move6a}
\end{figure}

We now want to show how this move~6 (and hence the light bulb trick)
can be performed when the $1$-handle has been replaced by a $(+1)$-surgery
along a Legendrian unknot $L_0$. Figure~\ref{figure:move6a}
shows that move~6 translates into a handle slide as in
Proposition~\ref{prop:bandsum}, where the connected sum $L\# L_0^+$
is not taken along a cusp, but along a pair of parallel strands;
cf.~\cite[Figure~2]{etho03} for this connected sum.
A connected sum of this kind with a Legendrian unknot can be undone
by a first Reidemeister move.

\begin{figure}[h]
\centering
\includegraphics[scale=0.45]{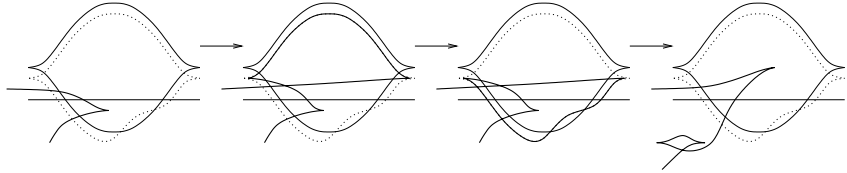}
  \caption{Move 6 via slide across meridional disc in surgered manifold.}
  \label{figure:move6b}
\end{figure}

The following alternative approach to this move, as illustrated
in Figure~\ref{figure:move6b}, may be worth noting. Here the construction
is based on the fact that $L_0^+$ is a Legendrian unknot
in the surgered manifold. Hence, we may assume that, in the
surgered manifold, $L_0^+$ bounds a disc foliated by Legendrian curves
as in Figure~\ref{figure:saucer-foliation}. This allows us to
move the upper strand of $L_0^+$ to the lower one by a
Legendrian isotopy across this new meridional disc
in the surgered manifold, as shown in the second step in
the sequence of Figure~\ref{figure:move6b};
the remaining isotopies are given by Reidemeister moves.

Both versions of handle addition used here for performing move~6
should be regarded as incarnations of the second Kirby move.

At first sight, it is a little surprising that if $L$ is a Legendrian
knot passing over a $1$-handle in a description of
some contact manifold $(M,\xi )$, then
there is a contactomorphism of $(M,\xi )$ sending $L$
to any of its stabilisations $S^k_+S^m_-L$. This can be seen from
Figure~\ref{figure:link-diagram-stab} by replacing $L_-$
with $S^k_+S^m_-L_-$, and $L_+$ with the Legendrian push-off
of $S^k_+S^m_-L_-$.

When we replace the $1$-handle in Figure~\ref{figure:link-diagram-stab}
by a contact $(+1)$-surgery along a Legendrian unknot~$L_0$, linked
once with $L$, this contactomorphism may be regarded as a Rolfsen twist
along~$L_0$ (see \cite[Figure~5.27]{gost99} for the
topological picture). Observe that the surgery coefficient of $L_0$
with respect to the surface framing of $L_0$ equals zero;
this coefficient remains zero under a Rolfsen twist.
\section{Long Legendrian knots}
We now return to the completion procedure $K\mapsto\hat{K}$ for
long Legendrian knots~$K$, where we regard $\hat{K}$ as a Legendrian knot
in $(S^3,\xist )$. The following theorem answers a
question of Fuchs and Tabachnikov~\cite[Section~3.5]{futa97}.

\begin{thm}
The completion map $K\mapsto\hat{K}$ sets up a one-to-one
correspondence between Legendrian isotopy classes of long Legendrian
knots in $(\R^3,\xist )$ and those of Legendrian knots in
$(S^3,\xist )$.
\end{thm}

\begin{proof}
The map $K\mapsto\hat{K}$ induces a surjective map $[K]\mapsto [\hat{K}]$
on Legendrian isotopy classes. We need to show that this last
map is injective.

Thus, let $K_0$ and $K_1$ be two long Legendrian knots in $(\R^3,\xist )$
whose completions $\hat{K}_0$, $\hat{K}_1$ are Legendrian isotopic
in $(S^3,\xist )$. This Legendrian isotopy extends to a
contact isotopy of $(S^3,\xist )$~\cite[Thm.~2.6.2]{geig08}.
Let $\phi_1$ be the time-$1$-map of that isotopy, which
we may assume to be the identity on an interval $I\subset
\hat{K}_0\cap \hat{K}_1$
containing the point $\infty :=S^3\setminus\R^3$. Our aim is to find
a contactomorphism $\phi_1'$ of $(\R^3,\xist )$ sending $K_0$ to $K_1$, and
with $\phi_1'=\mbox{\rm id}$ outside a compact set, say some large ball~$B$.
For then,
as at the end of the proof of Proposition~\ref{prop:trivial-long},
we see that $\phi_1'$ is in fact contact isotopic to
the identity relative to $\R^3\setminus B$. This translates into
a Legendrian isotopy between the long knots $K_0$ and $K_1$.

We provide two constructions of such a contactomorphism~$\phi_1'$.
The first one is based on the contact disc theorem; the second
construction, which is inspired by the proof of Lemma~5.3
in~\cite{etho03}, uses the classification of tight contact structures on solid
tori. We include the second argument because the method may have some
interest beyond the present application.

In new local coordinates $(x,y,z)$ around $\infty = (0,0,0)$
we may think of $I$ as the interval $(0,t,0)$, $|t|\leq\varepsilon$,
and $\xist$ as being given by $dz+x\, dy =0$. This will be the
set-up used for both arguments.

(1) Write $\phi_1(x,y,z)=(\tilde{x},\tilde{y},\tilde{z})$ and consider
the contact dilation
\[ \delta_s(x,y,z)=(sx,sy,s^2z).\]
Let $B_{\infty}$ be a small closed ball around $\infty =(0,0,0)$ with
$I\cap B_{\infty}$ contained in the interior of~$I$. As shown in
the proof of the contact disc theorem~\cite[Thm.~2.6.7]{geig08},
with $A:=\frac{\partial\tilde{z}}{\partial z}(0)\in\R^+$ and
$C:=\frac{\partial\tilde{y}}{\partial x}(0)$ we obtain a $1$-parameter
family of contact embeddings of $B_{\infty}$ by
\[ (x,y,z;s)\mapsto\left\{ \begin{array}{ll}
\delta_s^{-1}\circ\phi_1\circ\delta_s(x,y,z)&,\, s\in (0,1],\\[.5mm]
(Ax,Cx+y,Az-\frac{1}{2}ACx^2)&,\, s=0.\end{array}\right. \]
These contact embeddings of $B_{\infty}$ fix $I$ pointwise, and
the embedding for $s=0$ is contact isotopic to the identity
via contact embeddings fixing $I$ pointwise. The contact isotopy
extension theorem~\cite[Thm.~2.6.12]{geig08} then yields
the desired contactomorphism~$\phi_1'$.

(2) For $\delta >0$ sufficiently small, the embedded discs
$\phi_1(D_{y_0})$, $|y_0|\leq \varepsilon$, with
\[ D_{y_0}:=\{ (x,y_0,z)\co x^2+z^2\leq\delta^2\},\]
define a foliation transverse to $\partial_y$. So for each fixed pair
$(x_0,z_0)$, the holonomy map $y_0\mapsto y_0'$ sending $y_0$ to
the $y$-coordinate of the intersection of $\phi_1(D_{y_0})$ with the
$y$-line through $(x_0,0,z_0)$ defines, near $y_0=0$, a diffeomorphism onto
its image.

This allows us to apply the $1$-dimensional
disc theorem~\cite[Thm.~8.3.1]{hirs76} with
two parameters $(x_0,z_0)$ in order to isotope $\phi_1$ to
a diffeomorphism that sends each disc $D_{y_0}$ (possibly for some
smaller~$\delta$) into the plane $\R\times\{ y_0\}\times\R$,
fixing the origin. Next we apply the $2$-dimensional disc theorem
with one parameter $y_0$ in order to isotope this diffeomorphism
further to a diffeomorphism $\phi_2$ with $\phi_2 (K_0)=K_1$,
$\phi_2=\mbox{\rm id}$
near~$\infty$ and $\phi_2=\phi_1$ outside a slightly bigger neighbourhood.
Since $\phi_1$ was a contactomorphism, and the contact framing of
$I$ is given by $\partial_x$, we may in fact assume --- by the
neighbourhood theorem for Legendrian curves --- that $\phi_2$
is also a contactomorphism of a neighbourhood of $\hat{K}_0$ onto
a neighbourhood of~$\hat{K}_1$.

\begin{figure}[h]
\labellist
\small\hair 2pt
\pinlabel {Here $\phi_2=$ id} [t] at 47 5
\pinlabel {Here $T\phi_2(\xist )=\xist$} [l] at 247 0
\pinlabel $\phi_2$ [b] at 325 194
\pinlabel $\xist$ at 144 175
\pinlabel $T\phi_2(\xist)$ at 504 175
\pinlabel $T^2$ [bl] at 223 240
\pinlabel {$S^1\times D^2$} [br] at 59 214
\endlabellist
\centering
\includegraphics[scale=0.5]{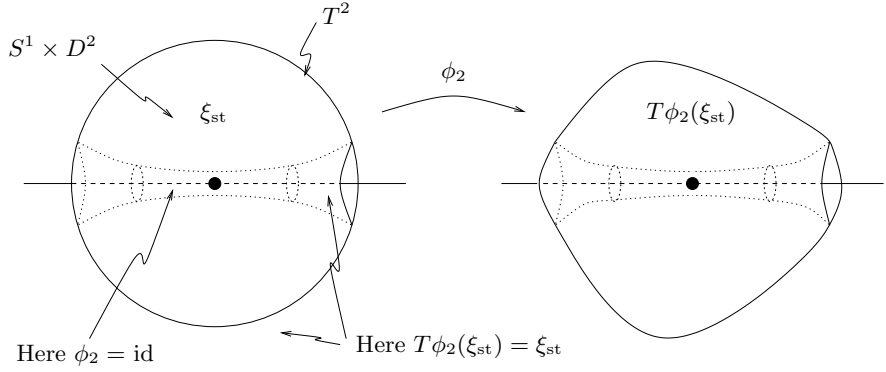}
  \caption{The diffeomorphism $\phi_2$ near $\infty$.}
  \label{figure:phi-torus}
\end{figure}

So we have a situation as illustrated
in Figure~\ref{figure:phi-torus}, with $\phi_2$ a contactomorphism
outside a solid torus $S^1\times D^2$. The boundary $T^2$ of that solid torus
is convex with characteristic foliation as shown in
Figure~\ref{figure:torus-boundary}; this has slope $-1$ (the dividing
curves are shown as dashed lines).

\begin{figure}[h]
\labellist
\small\hair 2pt
\pinlabel $\mu$ [t] at 217 -2
\pinlabel $\lambda$ [br] at -2 217
\endlabellist
\centering
\includegraphics[scale=0.35]{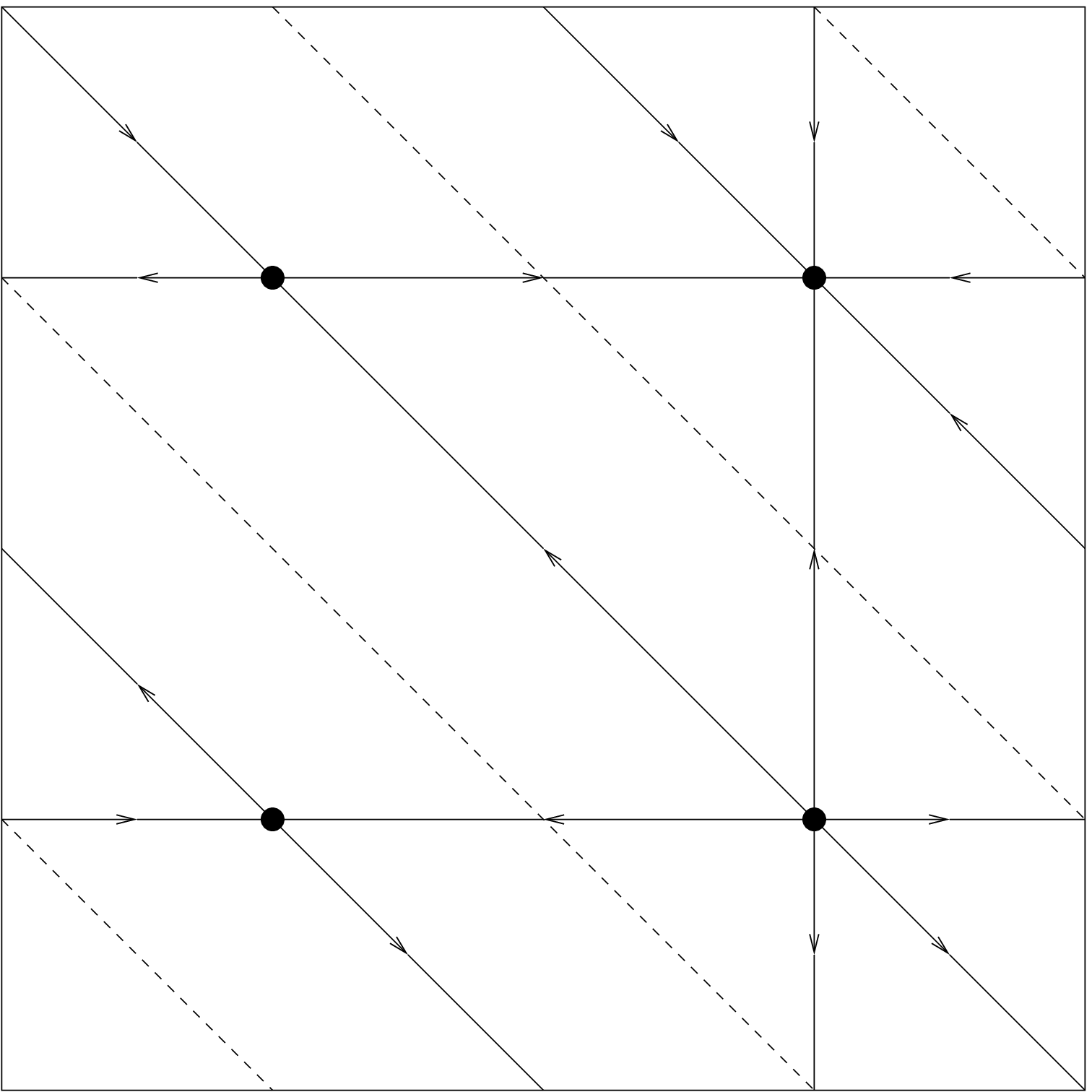}
  \caption{The characteristic foliation on $T^2$.}
  \label{figure:torus-boundary}
\end{figure}

Up to isotopy rel boundary, there is a unique tight contact structure
on the solid torus with this characteristic foliation on the
boundary~\cite[Prop.~4.3]{hond00}. So we find a contactomorphism
from $(\phi_2(S^1\times D^2),T\phi_2(\xist ))$ to
$(\phi_2(S^1\times D^2),\xist)$. Now extend this contactomorphism by
the identity outside $\phi_2(S^1\times D^2)$ to a
contactomorphism $\phi_3\co (S^3,T\phi_2(\xist))\rightarrow
(S^3,\xist)$. Then $\phi_1':=\phi_3\circ\phi_2$ has the
desired properties.
\end{proof}
\section{Destabilisation and the first Kirby move}
\label{section:destab}
A Kirby move of the first kind consists in adding to (or deleting from)
a given surgery diagram an unknotted circle, not linked with any component
of the given link, with framing $\pm 1$. This does not change the
topology of the $3$-manifold described by the surgery diagram.

\begin{figure}[h]
\labellist
\small\hair 2pt
\pinlabel $+1$ [bl] at 217 118
\pinlabel $+1$ [bl] at 577 118
\pinlabel $-1$ [bl] at 412 153
\endlabellist
\centering
\includegraphics[scale=0.3]{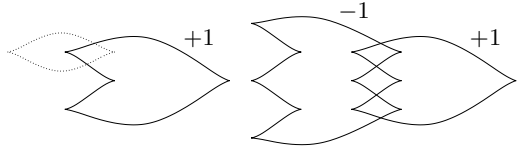}
  \caption{The `first Kirby move'.}
  \label{figure:kirby1}
\end{figure}

There is a not entirely satisfying analogue of this construction for
contact surgery diagrams
describing an overtwisted contact $3$-manifold. Here the move consists
in adding the link in Figure~\ref{figure:kirby1} to the given diagram.
This corresponds, as explained in~\cite{dgs04}, to forming the
connected sum with $S^3$ equipped with an overtwisted contact structure
that is homotopic, as a $2$-plane field, to~$\xist$; such a connected sum
does not change a given overtwisted contact manifold by
Eliashberg's theorem~\cite{elia89}, cf.~\cite[Section~4.7]{geig08}.
The dashed circle
indicates the boundary of an overtwisted disc, cf.~\cite[Fig.~2]{dgs04}.
(In the best of all possible worlds one should expect that it
were sufficient to add but the `shark' on the left,
corresponding to a single topologically $(-1)$-framed unknot. This, alas,
changes the Hopf $d_3$-invariant of the contact structure.)

\begin{figure}[h]
\labellist
\small\hair 2pt
\pinlabel $+1$ [tl] at 212 419
\pinlabel $+1$ [tl] at 502 419
\pinlabel $+1$ [tl] at 212 202
\pinlabel $+1$ [tl] at 502 202
\pinlabel $+1$ [tl] at 212 23
\pinlabel $+1$ [tl] at 502 23
\endlabellist
\centering
\includegraphics[scale=0.3]{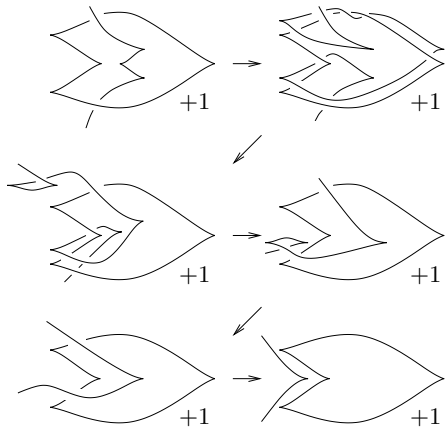}
  \caption{Destabilisation via Kirby moves.}
  \label{figure:kirby1-destab}
\end{figure}

Now suppose we are given a stabilised Legendrian knot in the original
diagram. When we perform the connected sum of such a knot with
the mentioned boundary of an overtwisted disc, we obtain (after a couple
of second Reidemeister moves) the situation
on the top left of Figure~\ref{figure:kirby1-destab}.
In the first row of that figure we perform a Kirby move of the
second kind; in the second and third row,
a move as in Figure~\ref{figure:meridian3}.
The remaining moves are the obvious Legendrian Reidemeister moves.
Thus, this figure gives an explicit demonstration of the well-known
result that the connected
sum of a Legendrian knot with the boundary of an overtwisted
disc (disjoint from the given knot) leads to a destabilisation,
cf.~\cite[pp.~316--319]{geig08}.

Another way of phrasing this result is that any Legendrian knot
with overtwisted complement (such `loose' knots will be
the topic of the next section) is Legendrian isotopic to
any of its stabilisations, with the appropriate number of
shark-like meridional curves added (as well as the curves
on the right of Figure~\ref{figure:kirby1}). This corresponds
to~\cite[Lemma~4.7]{dyma04}.

\begin{figure}[h]
\labellist
\small\hair 2pt
\pinlabel $S_+$ at 195 174
\pinlabel $S_-$ at 195 100
\pinlabel $+1$ at 465 164
\pinlabel $+1$ at 465 113
\endlabellist
\centering
\includegraphics[scale=0.3]{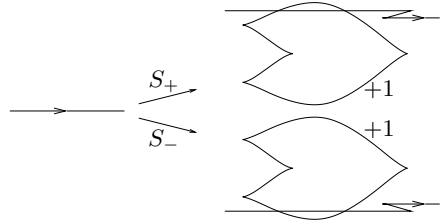}
  \caption{Stabilisations and respective position of shark.}
  \label{figure:shark-meridian}
\end{figure}

We leave it to the reader to verify the following simple observation. Both
positive and negative stabilisations can be effected by a
shark with his `mouth' on the left. One can always perform the
stabilisation along a strand oriented towards the right, say, by
first performing a first Legendrian Reidemeister move to create
an up or down cusp on the right. The sign of the stabilisation
affects the positioning of the shark-like meridian as shown
in Figure~\ref{figure:shark-meridian}.

\section{Loose knots}
A Legendrian knot $L$ in a contact $3$-manifold $(M,\xi )$ is called
{\em loose\/} if $\xi|_{M\setminus L}$ is overtwisted. The following theorem
is a generalisation of a result due to K.~Dymara~\cite{dyma04}.
When we speak of a Legendrian knot $L\subset (M,\xi )$ as a {\em framed\/}
knot, we always mean the framing defined by the contact structure.

\begin{thm}
\label{thm:loose}
Let $L_0,L_1$ be loose Legendrian knots in a (closed, connected)
contact $3$-manifold $(M,\xi )$, with an overtwisted disc
$\Delta\subset M\setminus (L_0\cup L_1)$. Assume that:
\begin{itemize}
\item[(i-a)] $L_0$ and $L_1$ are isotopic as framed knots;
\item[(i-b)] the topological knot type of $L_0,L_1$ contains infinitely many
non-isotopic framed knots;
\item[(ii-a)] $L_0$ and $L_1$ are Legendrian regularly homotopic as
Legendrian immersions;
\item[(ii-b)] the connected component of $L_0,L_1$ in the space of framed
immersions $S^1\rightarrow M$ contains infinitely many components
of the space of Legendrian immersions $S^1\rightarrow (M,\xi )$.
\end{itemize}
Then $L_0$ and $L_1$ are Legendrian isotopic.
\end{thm}

\begin{rem}
Condition (i-b) being satisfied for all knot types
in $M$ is equivalent to there being {\em no\/}
summands $S^1\times S^2$
in the prime decomposition of $M$, see the remark on p.~488
of~\cite{hopr90} and, for a proof, \cite[Thm.~2.0.5]{cher05}. One
direction of this equivalence is simple. The existence of such a summand is
equivalent to $M$ containing a non-separating $2$-sphere; the framing
of a knot $K$ intersecting such a $2$-sphere transversely in a single
point can be changed by the (topological) light-bulb trick by any
even number of twists, resulting in only two
different {\em framed\/} isotopy classes within the isotopy
class of~$K$.

Condition (ii-b) is equivalent to there being {\em no\/} class
$\alpha\in H_2(M)$ on which the Euler class of $\xi$ evaluates
non-trivially and that can be represented by a mapping $S^1\times S^1
\rightarrow M$ whose restriction to the second factor is a loop
freely homotopic to $L_0,L_1$, see~\cite[Prop.~3.1.4]{tche03}.
In particular, this condition is satisfied when
$\xi$ is trivial as a plane bundle (as was assumed in~\cite{dyma04}).
\end{rem}

\begin{proof}[Proof of Theorem~\ref{thm:loose}]
By \cite{dige04} we can represent $(M,\xi )$ by
a surgery diagram $\LL =\LL^-\sqcup\LL^+$ in $(S^3,\xist )$, i.e.\
a Legendrian link $\LL$ such that contact $(-1)$-surgery along the components
of $\LL^-$ and contact $(+1)$-surgery along the components of $\LL^+$
yields $(M,\xi )$. Thus, we may think of $L_0,L_1$ as
Legendrian knots in $(S^3,\xist )$ disjoint from~$\LL$.

By assumption, there is an isotopy from $L_0$ to $L_1$ (as framed knots)
in~$M$. This translates into a composition of isotopies in
$S^3\setminus\LL$ and handle slides over the $2$-handles
attached along~$\LL$. The latter can be performed as band-sums.

According to \cite[Thm.~4.4]{futa97}, any topological isotopy
between two Legendrian knots in $(S^3,\xist )$
can be converted into a Legendrian
isotopy of suitable stabilisations
of the two knots. The argument given there remains
valid inside $S^3\setminus\LL$. Similarly, any topological band-sum
can be performed as a Legendrian band-sum, provided we permit ourselves to
stabilise one of the summands (in our case: $L_0$)
as well as the Legendrian band, where stabilisations are necessary
for realising left-handed twists in the band.

Thus, we can Legendrian isotope a suitable stabilisation $S^k_+S^m_-L_0$
to some stabilisation $S^l_+S^n_-L_1$.
Since every stabilisation adds a negative twist to the contact framing,
conditions (i-a) and (i-b) imply that $k+m=l+n$. Similarly,
conditions (ii-a) and (ii-b) imply that $k-m=l-n$. This follows
from the fact that the Legendrian regular homotopy class
of a parametrised closed curve $\gamma\co S^1\rightarrow (M,\xi )$
is determined by the homotopy class of the bundle map
$T\gamma\otimes\C\co TS^1\otimes\C\rightarrow\xi$, cf.\
\cite[Section~6.3]{geig08}, and the effect of a positive or negative
stabilisation on this bundle map which, being a local modification,
can be understood in the same way as the effect of a stabilisation on the
rotation number. We conclude that $k=l$ and $m=n$.

The Legendrian isotopy between $S^k_+S^m_-L_0$ and $S^k_+S^m_-L_1$
extends to a contact isotopy $\phi_t$, $t\in [0,1]$, of $(M,\xi )$.
Since $\phi_t(L_0)$, $t\in [0,1]$, defines a Legendrian isotopy from
$L_0$ to $\phi_1(L_0)$, we may replace $L_0$ in the sequel of the argument
by~$\phi_1(L_0)$. (Observe that we may construct $\phi_t$ in such a way that
it remains stationary on the overtwisted disc~$\Delta$, so the complement
of $L_1$ and this new $L_0$ still contains~$\Delta$.)

Thus, we may assume that $S^k_+S^m_-L_0$ and $S^k_+S^m_-L_1$ actually
coincide for some specific choice of these stabilisations.
At any rate, we also
find some segment where $L_0$ and $L_1$ coincide.
Since there is an overtwisted disc in the complement of $L_0\cup
L_1$, we may perform first Kirby moves in $M\setminus (L_0\cup L_1)$
without changing the situation (up to contactomorphism). Thus,
by the discussion in Section~\ref{section:destab}, the knots
$L_0,L_1$ are Legendrian isotopic in $(M,\xi )$
to their stabilisations $S^k_+S^m_-L_0=S^k_+S^m_-L_1$ with the
appropriate number of shark-like meridians added; those
meridians can be added in a neighbourhood of the segment
where $L_0$ and $L_1$ coincide and thus may be taken to
be identical for $L_0$ and~$L_1$. This proves
that $L_0$ and $L_1$ are Legendrian isotopic.
\end{proof}

\begin{rem}
Conditions (i-a) and (ii-a) in Theorem~\ref{thm:loose} are obviously
necessary. The examples in Section~4 of \cite{tche03} show
that (i-b) and (ii-b) are likewise necessary.
\end{rem}

With the help of the contact surgery presentation theorem from
\cite{dige04} one can also formulate a result about isotopies
of Legendrian knots $L_0,L_1$ in an arbitrary contact $3$-manifold
$(M,\xi )$. Think of $(M,\xi )$ as being represented
by a surgery diagram $\LL =\LL^-\sqcup\LL^+$ in $(S^3,\xist )$,
and of $L_0,L_1$ as knots in $S^3\setminus\LL$. We visualise $L_0,L_1$
and the link $\LL$ in the front projection. Then $L_0$ and $L_1$
are Legendrian isotopic in $(M,\xi )$ if and only if the front of
$L_0$ can be turned into the front of $L_1$ by
Legendrian Reidemeister moves of the first three
kinds (in the complement of~$\LL$) and by second Kirby moves.
This follows in analogy with the topological case and
from the corresponding result for $\LL =\emptyset$~\cite{swia92}.
Any Legendrian isotopy
that does not pass through the belt sphere of one of the $2$-handles
corresponding to $\LL$ can be made disjoint from those
$2$-handles and thus takes place in $S^3\setminus\LL$; a Legendrian
isotopy with a single transverse passage through one of the belt spheres
corresponds to a second Kirby move.
\begin{ack}
Some of this research was carried out during a stay of F.~D. at the
Mathematical Institute of the Universit\"at zu K\"oln, supported by
grant no.\ 10631060 of the National Natural Science Foundation 
of China. H.~G.\ is partially supported by DFG grant
GE 1245/1-2 within the framework of the Schwer\-punkt\-pro\-gramm 1154
``Globale Differentialgeometrie''.
\end{ack}

\end{document}